\documentclass[journal]{IEEEtran}
\usepackage{graphics, graphicx, subfigure}
\usepackage{amsmath,amssymb}
\usepackage{amsthm}
\usepackage{bm} 
\theoremstyle{plain}
\newtheorem{theorem}{Theorem}[section]

\theoremstyle{definition}

\newtheorem{definition}[theorem]{Definition}

\newtheorem{problem}[theorem]{Problem}

\begin{document}

\title{Approximation by Spline Curves: \\ towards an Application to Cognitive Neuroscience}

\author{Maria-Laura~Torrente,
        Stefano~Anzellotti, Chiara Finocchiaro
        and~Claudio~Fontanari
\thanks{This research project was partially supported by GNSAGA of INdAM, by PRIN 2010--2011
"Geometria delle variet\`a algebriche", and by FIRB 2012 "Moduli spaces and Applications".}
\thanks{M. Torrente is with the Dipartimento di Matematica, Universit\`a di Genova, 16146 Genova, Italy
(e-mail: torrente@dima.unige.it).}
\thanks{S. Anzelotti is with the Department of Brain and Cognitive Sciences, 
Massachussetts Institute of Technology, Cambridge, MA 02139 USA
(e-mail: anzellot@mit.edu).}
\thanks{C. Finocchiaro is with the Dipartimento di Psicologia e Scienze Cognitive,
 Universit\`a degli Studi di Trento, 38068 Rovereto, Italy
(e-mail: chiara.finocchiaro@unitn.it).}
\thanks{C. Fontanari is with the Dipartimento di Matematica and
Centro Interdipartimentale Mente/Cervello-CIMeC,
Universit\`a degli Studi di Trento, 38123 Trento, Italy
(e-mail: fontanar@science.unitn.it).}}


\maketitle

\begin{abstract}
We present a procedure to approximate a plane contour by piecewise polynomial functions, depending 
on various parameters, such as degree, number of local patches, selection of knots.
This procedure aims to be adopted to study how information about shape is represented.
\end{abstract}


\IEEEpeerreviewmaketitle

\section{Introduction}

The shape of an object is a fundamental source of information for recognition. While color and texture information are also important 
(see \cite{Yip}), line drawings, which discard texture and color,  are usually sufficient to recognize an object. 
Furthermore, infants generalize object labels on the basis of shape similarity rather than color similarity (see \cite{Graham}), 
suggesting that shape plays a key role in categorization. 
In order to understand how we recognize and categorize objects, therefore, it is crucial to study how information about shape is
represented. 

Since the space of shapes is infinite dimensional, the brain likely approximates it using a lower number of dimensions
 in order to make the problem more tractable.
For comparison, it can be useful to consider the example of color. Like the space of shapes, the space of light spectra is infinite 
dimensional. However, with three dimensions it is possible to model accurately the space of perceptually discernible colors 
(see \cite{Krauskopf}). In the same vein, we can look for finite dimensional spaces that encompass the perceptually discernible shapes.

Recent work (see \cite{Hong}) found that a model with approximately $47$ dimensions explains accurately a large amount 
of the variance in the object recognition errors and similarity judgments of human participants. However, models of this type
leave open the question of whether the dimensionality of shape space is the same in the neighborhoods of all shapes. 
In addition, the most relevant dimensions for recognition  might differ in the neighborhood of different shapes. 
For example, the relevant dimensions for discriminating between individual faces and for discriminating between buildings 
might not be the same, and thus different low-dimensional approximations of the space of shapes could be used in the two 
neighborhoods. This possibility is compatible with the differential involvement of different brain regions within temporal cortex 
in the processing of different object categories (see \cite{Reddy}). 
Our interdisciplinary research project aims to introduce techniques to individuate lower dimensional subspaces that 
locally approximate the space of perceptually discernible shapes.

A couple of crucial remarks are in order here. First, our approach aims to model shape discriminability rather than object recognition accuracy. In other words, it is a model of perception rather than a model of perceptual categorization (in this respect, it is quite similar to the case of color space). The choice of investigating discriminability is based on the idea that 
a common representation of shape underlies both shape discrimination and judgements of perceptual similarity that are 
at the basis of perceptual categorization. 
This idea is supported by the evidence currently available (see \cite{Hong}). 
Next, our goal is to represent  shapes locally rather than constructing a single 
lower dimensional space underlying the representation of all shapes. This choice is a consequence of the aim to model perceptual discriminability: different local spaces seem necessary to account for the large variety of shapes we can perceptually discriminate.

For modeling shape space in the neighborhood of a fixed shape spline curves will be used.
Spline curves can approximate contours closely, they enjoy elegant mathematical properties, and 
their computational complexity is relatively low. For all these reasons, splines are widely considered
to be an ideal tool to approximate signals (see for instance \cite{Unser}). In particular, the use 
of B-splines for digital signal processing has a long history (at least since \cite{Hou}). 
In this paper we employ spline curves to address the following approximation problem: given a 
plane contour $\mathcal C$ construct a (parametric) spline curve~$C$ of a fixed degree $d$ 
well approximating $\mathcal C$.
The procedure (described in Section \ref{algorithmSec} and based on Section \ref{basicSpline}) 
is divided into three parts: use of strengthened 
edge detection techniques to extract the points of the contour $\mathcal C$; computation of the spline curve $C$, 
performed by solving a Least Squares Approximation problem (see Problem \ref{LSApproxProblem}
and subsequent discussion); plot of the computed spline curve (and, eventually, comparison with~$\mathcal C$).
The spline approximation part (second part) is the core of the presented procedure and, as it is well known 
to the experts (see for instance \cite{He}), an important issue for it is the choice of knots.
At the end of Section \ref{basicSpline} we recall three of the most used parametrization methods 
(uniform, chord length and centripetal). The chord length parametrization method is then chosen to show 
a concrete example (see figures \ref{origImage}-\ref{caso4}).
 

\section{Basics of Spline Approximation}\label{basicSpline}
The main reference for this section is \cite{Lyche}. 

B\'ezier curves and, after their introduction in 1946 (see \cite{Schoenberg}), 
spline curves have been widely used to construct smooth curves from 
a given set of points in an efficient and numerically stable way.
Their geometrical construction is essentially based on recursive convex 
combinations of curves of smaller degrees. Namely:

\begin{definition}\label{splinecurve}
Let $n,d$ be positive integers, with $n \ge d+1$.
Let $p_1,\ldots,p_n \in \mathbb R^2$ be $n$ control points,
and let $\tau = (t_i)_{i=2}^{n+d}$ be the knot vector, which is assumed 
to be a nondecreasing sequence of $n+d-1$ real numbers, that is, 
$t_2 \le \ldots \le t_{n+d}$.
The functions ${\bf p}_{i,j}(t)$ are recursively defined in the following 
way: we set ${\bf p}_{i,0}(t) = p_i$, for $i=1,\ldots, n$ and 
\begin{eqnarray*}
{\bf p}_{i,j}(t) = \frac{t_{i+d-j+1}-t}{t_{i+d-j+1}-t_i} {\bf p}_{i-1,j-1}(t) 
 +\frac{t-t_i}{t_{i+d-j+1}-t_i} {\bf p}_{i,j-1}(t) 
\end{eqnarray*}
for $j=1,\ldots,d$ and $i=j+1,\ldots,n$, 
where possible division by zero is resolved by the convention 
that a division by zero is zero.
The {\it (parametric) spline curve~$f(t)$ of degree $d$ with control points 
$(p_i)_{i=1}^n$ and knots $\tau=(t_i)_{i=2}^{n+d}$} is defined as
$$
f(t)= \sum_{i=d+1}^n {\bf p}_{i,d}(t) B_{i,0}(t)
$$ 
where $B_{i,0}(t)$ are the piecewise constant functions defined  by 
\begin{eqnarray*}
 B_{i,0}(t) = \left \{ 
 \begin{array}{lll}
 1 & & \textrm{if } t_i \le t < t_{i+1}\\
 0 & & \textrm{otherwise} 
\end{array} \right.
\end{eqnarray*}
\end{definition} 

Due to their special construction, spline curves enjoy many important 
geometrical properties, such as local dependence on the control points, 
being contained in the control polygon, and so on.

A crucial fact is that every spline curve $f(t)$ can be written 
as a linear combination of special univariate polynomials (the B-splines,
see Definition~\ref{B-splines}) with the control points as coefficients
(see Theorem~\ref{splinesAsBsplines}).
\begin{definition}\label{B-splines}
Let $m,d$ be nonnegative integers, with $m \ge d+2$, and 
let $\tau = (t_i)_{i=1}^m$ be a nondecreasing sequence of real numbers.
Let $i=1,\ldots,m+d+1$; the {\it $i$th B-spline} of degree $d$ 
with knots $\tau(t_i)_{i=1}^m$, denoted by $B_{i,d}(t)$, is recursively defined 
in the following way: for $i=1,\ldots, m-1$, we set 
\begin{eqnarray*}
B_{i,0}(t)=\left \{ \begin{array}{lll}
1 & & \textrm{if } t_i \le t \le t_{i+1}\\
0 & & \textrm{otherwise}
\end{array} \right.
\end{eqnarray*}
Then, for each $j=1,\ldots,d$ and $i=1, \ldots,m-j-1$
we define
\begin{eqnarray*}
B_{i,j}(t) = \frac{t-t_i}{t_{i+j}-t_i}B_{i,j-1}(t) 
 +\frac{t_{i+1+j}-t}{t_{i+1+j}-t_{i+1}} B_{i+1,j-1}(t) 
\end{eqnarray*}
where, as above, possible division by zero is resolved by the convention 
that a division by zero is zero.
\end{definition}

\begin{theorem}\label{splinesAsBsplines}
Let $n,d$ be positive integers, with $n \ge d+1$.
Let $p_1,\ldots,p_n \in \mathbb R^2$ be $n$ control points,
and let $\tau = (t_i)_{i=1}^{n+d+1}$ be a nondecreasing sequence of real
numbers, that is, $t_1 \le \ldots \le t_{n+d+1}$.
The (parametric) spline curve $f(t)$ of degree $d$ with control points 
$(p_i)_{i=1}^n$ and knots $\tau= (t_i)_{i=1}^{n+d+1}$ can be written as:
$$
f(t)=\sum_{i=1}^n p_i B_{i,d}(t)
$$
where $B_{i,d}(t)$ is the $i$th B-spline of degree $d$ 
with knots $\tau$.
\end{theorem}

The notion of linear combination of B-splines, also called spline function,
is formalized in the following definition.

\begin{definition}\label{splineFunctions}
Let $n,d$ be positive integers, with $n \ge d+1$.
Let $\tau = (t_i)_{i=1}^{n+d+1}$ be a nondecreasing sequence
of reals, and let $B_{1,d}(t), \ldots, B_{n,d}(t)$ be the $n$ B-splines
of degree $d$ with knots $\tau$.
The linear space of all linear combinations of 
$B_{1,d}(t), \ldots, B_{n,d}(t)$ is called the {\it spline space
of degree $d$ with knots $\tau$} and denoted by:
\begin{eqnarray*}
\mathbb S_{d, \tau} = \left\{\sum_{i=1}^n c_i B_{i,d} \: \big|\: c_i \in \mathbb R,
i=1,\ldots,n \right\}
\end{eqnarray*}
An element $f=\sum_{i=1}^n c_i B_{i,d}$ of $\mathbb S_{d, \tau}$ is 
called a {\it spline function of degree $d$ with knots $\tau$}, and 
$(c_i)_{i=1}^n$ are called the {\it B-spline coefficients} of $f$.
\end{definition}

Now we consider the approximation problem.
The method is based on the standard least squares approach, 
in which classically the function to minimize is the sum 
of squared errors committed at the given data points.  
The problem can be formulated as follows.

\begin{problem}\label{LSApproxProblem}
{\bf (Least Squares Approximation Problem)}
Let $m > n$; given data $(x_k, y_k)_{k=1}^m$ and 
a spline space $\mathbb S_{d,\tau}$ whose knot vector
$\tau=(t_i)_{i=1}^{n+d+1}$ satisfies $t_{i+d+1}>t_i$ for $i=1,\ldots,n$, 
find a spline function $f=\sum_{i=1}^n c_i B_{i,d}(t) \in \mathbb S_{d,\tau}$
which solves the minimization problem 
\begin{eqnarray*}
\sum_{k=1}^m (y_k - f(x_k))^2 &=&
\sum_{k=1}^m \left(y_k - \sum_{i=1}^n c_i B_{i,d}(x_i) \right)^2 = \\
&=& \min_{g \in \mathbb S_{d,\tau}} \sum_{k=1}^m (y_k - g(x_k))^2
\end{eqnarray*}
\end{problem}

It is easy to see that this is a LS problem which 
can be expressed in matrix form as follows:
find ${\bm c} =(c_1,\ldots,c_n)^t$ s.t.
\begin{eqnarray}\label{approxWithNorms}
\|{\bm y} - B{\bm c}\|^2_2 = 
\min_{{\bm \alpha}=(\alpha_1,\ldots,\alpha_n) \in \mathbb R^n} \|{\bm y} - B {\bm \alpha} \|_2^2
\end{eqnarray}
where
\begin{eqnarray*}
B= \left (\begin{array}{ccc}
B_{1,d}(x_1) & \ldots & B_{n,d}(x_1)\\
\vdots & \ddots & \vdots\\
B_{1,d}(x_m) & \ldots & B_{n,d}(x_m)
\end{array} \right) \in \textrm{Mat}_{m \times n}(\mathbb R)
\end{eqnarray*}
is the coefficient matrix, also called {\it B-spline collocation
matrix}, and
\begin{eqnarray*}
{\bm c}= \left (\begin{array}{c}
c_1 \\ \vdots \\ c_n
\end{array}\right) \in \textrm{Mat}_{n \times 1}(\mathbb R), 
 \;
{\bm y}= \left (\begin{array}{c}
y_1 \\ \vdots \\ y_m
\end{array}\right) \in \textrm{Mat}_{m \times 1}(\mathbb R) 
\end{eqnarray*}
It is well-known that solving (\ref{approxWithNorms})
is equivalent to solve the normal system:
\begin{eqnarray}\label{normalSys}
B^tB {\bm c}=B^t {\bm y}
\end{eqnarray}

Note that Problem \ref{LSApproxProblem}
has always a solution. Moreover, there are explicit conditions under which 
the matrix $B^tB$ is invertible, and consequently system (\ref{normalSys}), 
as well as Problem \ref{LSApproxProblem}, has exactly one solution
(see for instance \cite{Lyche}, Theorem 5.23).

Finally, we address the following problem: suppose 
we are given a set $\mathbb X =\{p_1,\ldots,p_m\}$ of $m$ points 
in $\mathbb R^2$, and we want to construct a (parametric)
spline curve $\mathcal C$ of degree $d$ that approximates the points.
From Problem \ref{LSApproxProblem}, recalling Definition \ref{splineFunctions}, 
in order to define a spline space $\mathbb S_{d, \tau}$ 
it is clear that a knot vector $\tau$ is required. Obviously, different choices 
of the knot vector lead to different spline spaces and, consequently, 
to different approximating spline curves. 
In the literature (see for instance \cite{Lee}), there are various
competing parametrization methods, such as:
\begin{enumerate}
\item {\it Uniform:} $t_1=0$ and $t_i=\frac{i}{m}$, for $i=2,\ldots,m$.
\item {\it Chord length:} $t_1=0$ and 
$t_i=t_{i-1} + \frac{\|p_i- p_{i-1}\|_2}{\sum_{j=1}^n \|p_{j+1}- p_j\|_2}$, 
for $i=2,\ldots,m$.
\item {\it Centripetal:} $t_1=0$ and 
$t_i=t_{i-1} + \frac{\|p_i- p_{i-1}\|^{1/2}_2}{\sum_{j=1}^n \|p_{j+1}- p_j\|^{1/2}_2}$, 
for $i=2,\ldots,m$.
\end{enumerate}
In this paper we mainly focused on the chord length method, 
though the procedure for approximating plane contours using 
spline curves (described in Section \ref{algorithmSec}) 
has been implemented with all three parametrization methods.

\section{Algorithm Description and Implementation}\label{algorithmSec}
In this section we describe an algorithmic procedure for approximating 
a given plane contour  using spline curves. 
The procedure can be divided into three parts:
\begin{enumerate}
\item[I.] acquisition of a digital image reproducing a silhouette; extraction 
of the contour $\mathcal C$ by strengthened edge detection 
techniques;
\item[II.] computation of the spline curve, approximation of the contour $\mathcal C$;
\item[III.] plot of the computed spline curve (and eventually, comparison with $\mathcal C$).
\end{enumerate}
Accordingly, the implementation of the algorithm is divided into three steps:
the first and the third step, corresponding to parts I and III, are 
implemented in MatLab, using well-established built-in functions; the 
second step, corresponding to part II, is implemented in C++ language, 
using routines both of CoCoALib  \cite{CoCoALib}, a GPL C++ library -- the 
mathematical kernel for the computer algebra system CoCoA-5, and 
of the numerical library GSL - GNU Scientific Library \cite{GSL}.
Note that part~II is the core of the algorithm: based on Section \ref{basicSpline}, 
and in particular on Problem \ref{LSApproxProblem} and the subsequent
discussion, its implementation requires as input:
\begin{itemize}
\item[I1.] an affine set of points in $\mathbb R^2$
(resulting from the edge detection on the input image);

\item[I2.] a parameter $n$, corresponding to the number of polynomial patches;

\item[I3.] a parameter $d$, corresponding to the degree of polynomials;

\item[I4.] a parametrization method (see Section \ref{basicSpline}).
\end{itemize}
It returns as output:
\begin{itemize}
\item[O1.] a coefficient vector defining the piecewise polynomial function made up 
of $n$ patches of degree $d$;

\item [O2.] the corresponding Least Square Error.
\end{itemize}
 
\begin{figure}[htbp]
\centering%
\subfigure[Original image.\label{origImage}]%
          {\includegraphics[width=0.24\textwidth]{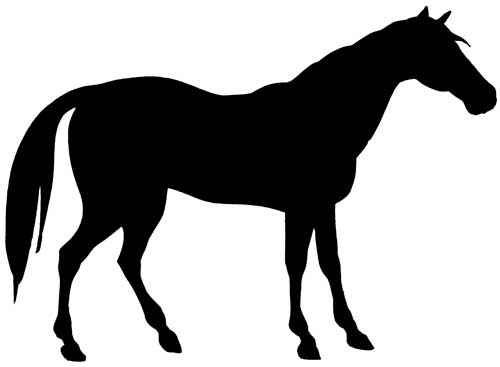}}
\subfigure[Points of the edge detection.\label{edge}]%
          {\includegraphics[width=0.24\textwidth]{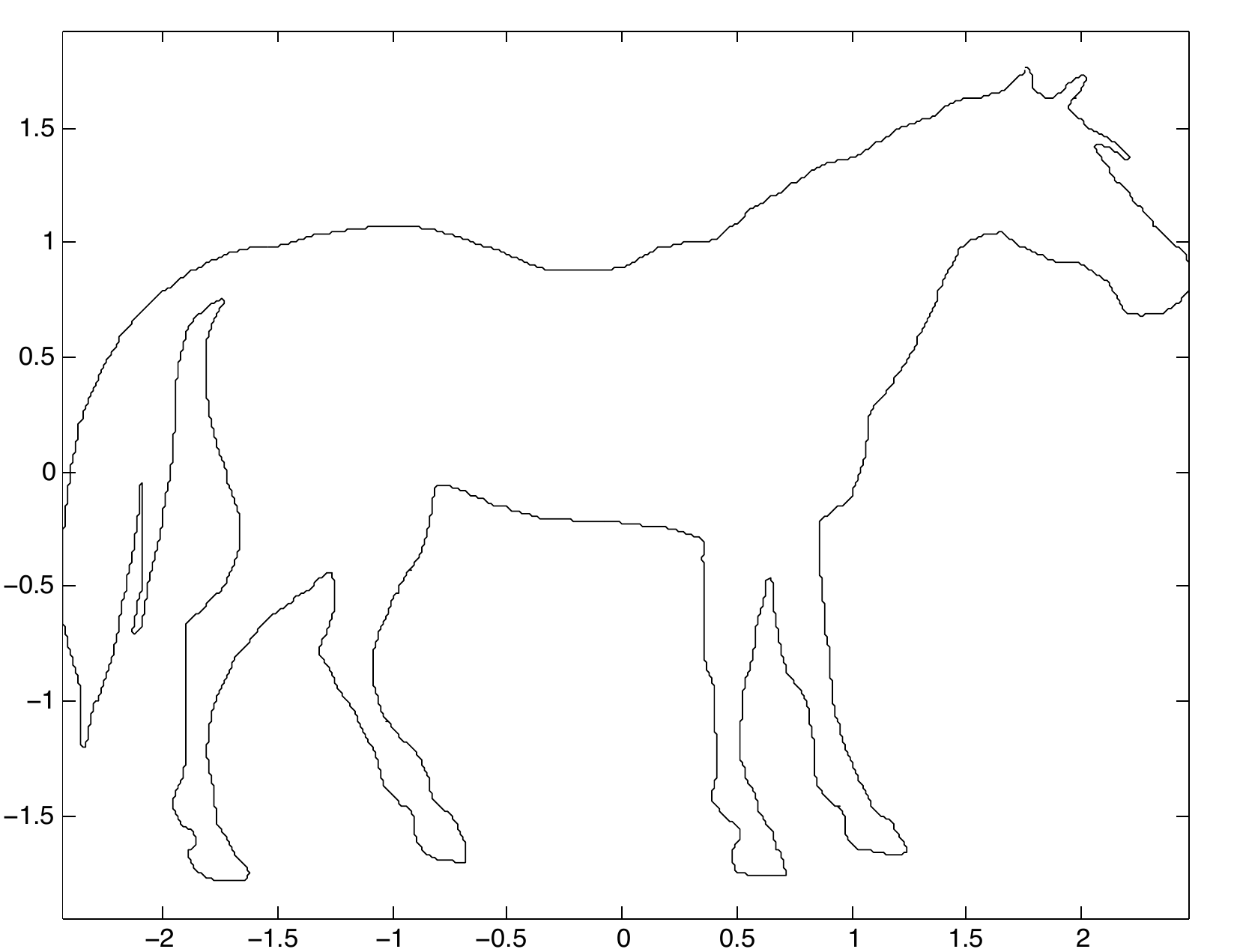}}
\subfigure[Values $n=100$, $d=2$.\label{caso1}]%
          {\includegraphics[width=0.24\textwidth]{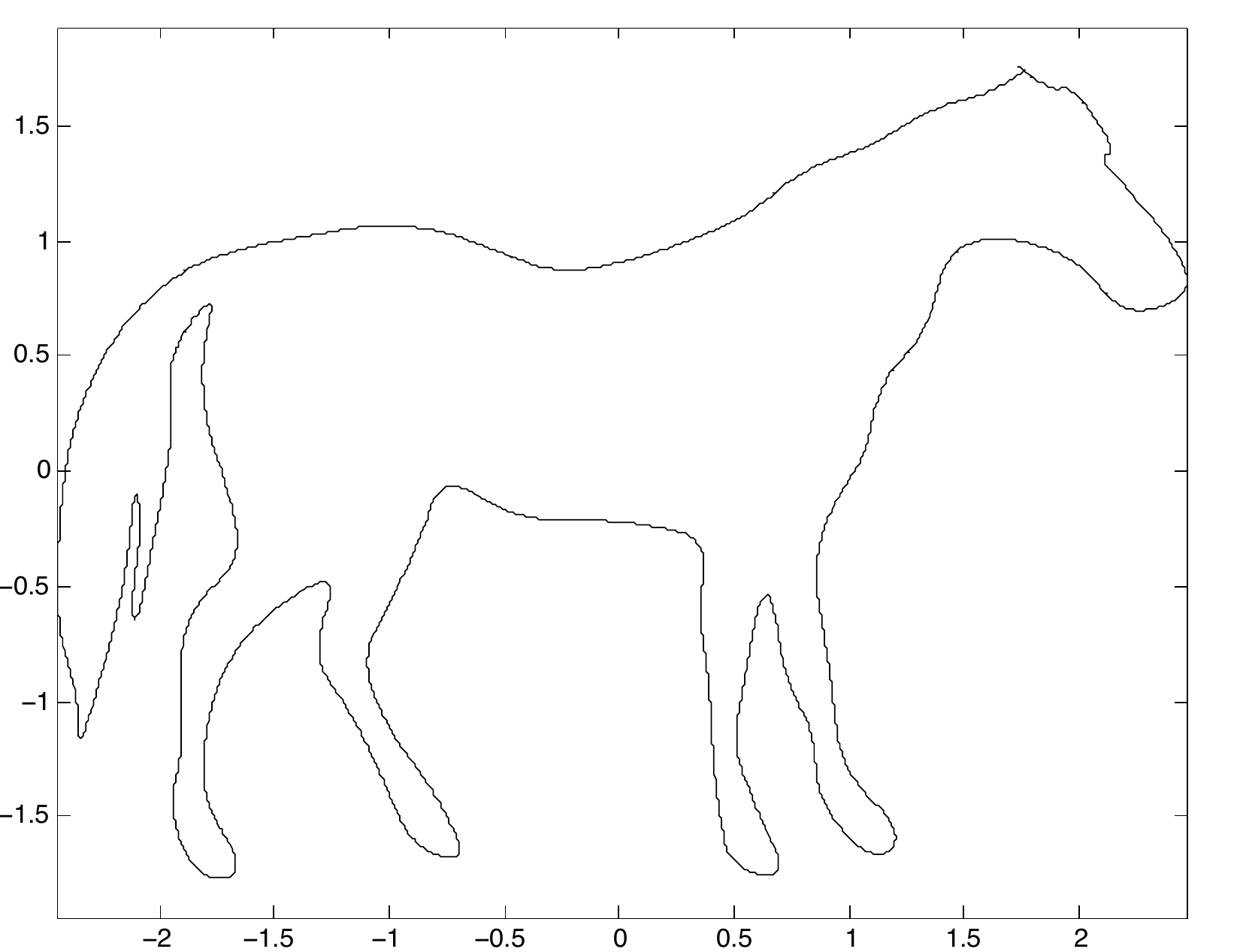}}
\subfigure[Values $n=200$, $d=2$.\label{caso2}]%
          {\includegraphics[width=0.24\textwidth]{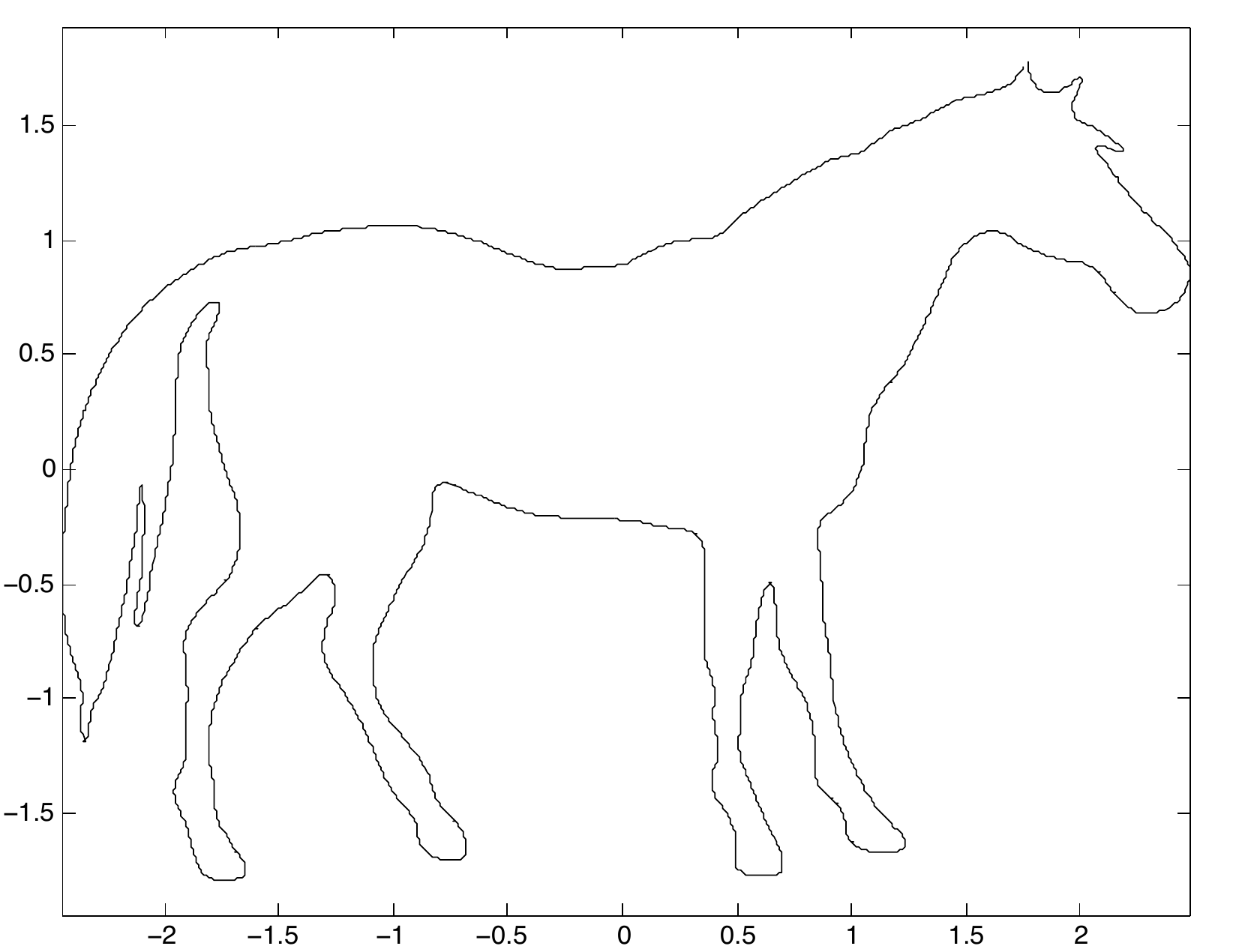}}
\subfigure[Values $n=100$, $d=3$.\label{caso3}]%
          {\includegraphics[width=0.24\textwidth]{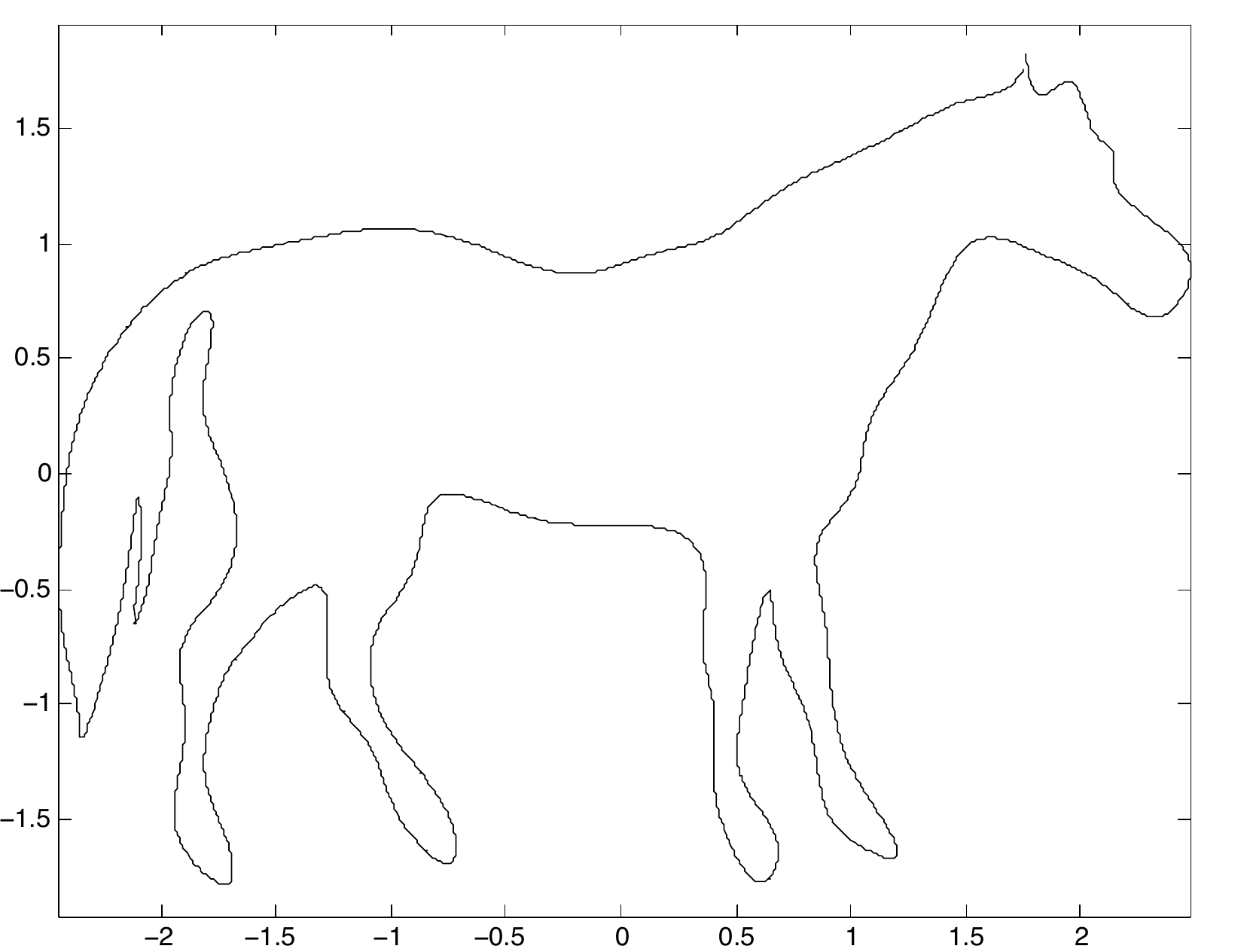}}
\subfigure[Values $n=200$, $d=3$.\label{caso4}]%
          {\includegraphics[width=0.24\textwidth]{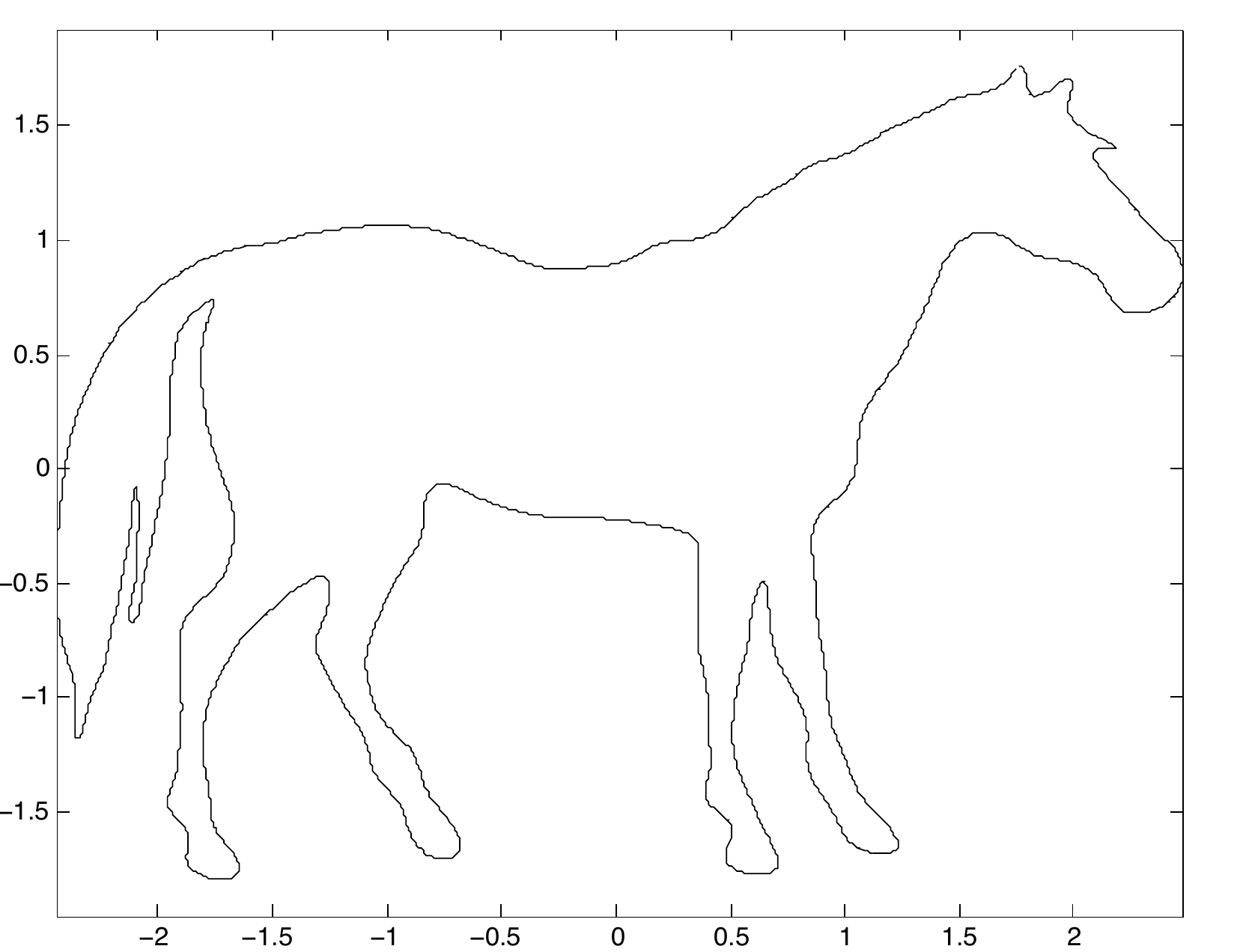}}
\centering\caption{Silhouette of a horse}
\end{figure}

%
%
%
Just to give an idea, we consider the silhouette of a horse, represented 
in Figure \ref{origImage}. Using classical edge detection techniques, 
the contour of the horse is extracted (see Figure \ref{edge}) and stored as
a set of $2713$ points in $\mathbb R^2$.
The chord length parametrization method is chosen (see Section \ref{basicSpline}), 
and contour approximations are obtained through computations 
with different values of the parameters, namely
$n=100, d=2$ (see Figure \ref{caso1}), $n=200, d=2$ (see Figure \ref{caso2}),
$n=100, d=3$ (see Figure \ref{caso3}) and $n=200, d=3$ (see Figure \ref{caso4}). 
We observe that there is a substantial difference in the approximation quality 
comparing the cases $n=100$ and $n=200$, whereas we notice a slight
improvement passing from the cases with $d=2$ to $d=3$. 

%
%

\section{Conclusion}
We have presented a novel procedure to approximate a plane contour by piecewise polynomial functions, depending 
on various parameters (degree, number of local patches, selection of knots). In order to optimize our choice of parameters, we plan to perform a series of behavioural experiments by adopting the following procedure. First, a set of contours will 
be selected, representing both shapes of daily life objects and meaningless curves. Next, every contour will be approximated by different splines, obtained by varying parameters in such a way that each approximation is perceptually closed to the original shape. Then, participants will be asked to perform a same/different task. Finally, the experimental results shall allow us to individuate a lower dimensional spline space approximating the space of perceptually discriminable shapes in a neighborhood of a fixed shape. At this point, it will be possible to ask questions about shape representations such as: 
\begin{enumerate}
\item Is the dimensionality of the space of perceptually discriminable shapes the same for neighborhoods of different shapes? 
\item In particular, does the dimensionality of the space of perceptually discriminable shapes change as a function of experience? 
\item As a consequence, do neighborhoods of real objects have higher dimensionality than neighborhoods of meaningless curves?
\end{enumerate}

\end{document}